\documentclass[12pt]{article}
\usepackage{amsmath,amssymb}
\usepackage{amsfonts,graphicx,latexsym}
\usepackage[notref,notcite]{} 
\usepackage{color,graphicx}        % standard LaTeX graphics tool

\numberwithin{equation}{section}

\newcommand{\adbc}{}
\newenvironment{proofsect}[1] 
{\vskip0.1cm\noindent{\bf #1.}\hskip0.5cm} 
 
\newcommand{\halmos}{\rule{1ex}{1.4ex}}
\newcommand{\eqa}{\begin{eqnarray}}
\newcommand{\ena}{\end{eqnarray}}
\newcommand{\eq}{\begin{equation}}
\newcommand{\en}{\end{equation}}
\newcommand{\eqs}{\begin{eqnarray*}}
\newcommand{\ens}{\end{eqnarray*}}
\def\qed{$\halmos$} 
\def\ep{\hfill\qed}
\def\Ref#1{{\rm (\ref{#1})}}
 
\def\ti{\to\infty} 
\def\a{\alpha}

\def\G{\Gamma} 
\def\de{\delta} 
\def\e{\varepsilon} 
\def\l{\lambda} 

\def\r{\varrho} 
\def\S{\Sigma}

\def\o{\omega}

\def\L{\Lambda} 
\newcommand{\R}     {\mathbb{R}} 
\newcommand{\Z}     {\mathbb{Z}}

\newcommand{\E}     {\mathbb{E}} 
 \newcommand{\floor}[1]{\left\lfloor #1 \right\rfloor}

\def\1{{\mathchoice {1\mskip-4mu\mathrm l}      % Blackboard bold 1 
{1\mskip-4mu\mathrm l} 
{1\mskip-4.5mu\mathrm l} {1\mskip-5mu\mathrm l}}} 
\newcommand{\ssup}[1] {{{\scriptscriptstyle{({#1}})}}} 
\def\comment#1{} 
 
\newtheorem{theorem}{Theorem}[section] 
 
\newtheorem{proposition}[theorem]{Proposition} 
\newtheorem{corollary}[theorem]{Corollary} 
\newtheorem{remark}[theorem]{Remark}

\newtheorem{example}[theorem]{Example}

\newcommand{\heap}[2]{\genfrac{}{}{0pt}{}{#1}{#2}} 
 
\newcommand{\s}{\sigma}
\renewcommand{\d}{{\rm d}} 
 
\newcommand{\eps}{\varepsilon}

\newcommand{\Bcal}  {{\mathcal B}}
 
\newcommand{\Dcal}   {{\mathcal D }} 
 
\newcommand{\Fcal}   {{\mathcal F }} 
\newcommand{\Gcal}   {{\mathcal G }}

\renewcommand{\e}   {{\operatorname e }}

\def\nin{\noindent}
\def\ignore#1{}
\def\Eq{\ =\ }
\def\Le{\ \le\ }
\def\Def{\ :=\ }
\def\e{\eps}
\def\iid{i.i.d.}
\def\bbP{\mathbb{P}}
\def\bbE{\mathbb{E}}
\def\nat{\mathbb{N}}
\def\bP{\mathbf{P}}
\def\bX{\mathbf{X}}
\def\bA{\mathbf{A}}

\def\giv{\,|\,}
\def\up{^{(\s)}}
\def\nti{\floor{nt}}
\def\de{\delta}
\def\parent#1{#1^{-1}}
\def\mmu{\beta}
\def\ki{{1}}
\def\kt{{2}}
\def\oke{\overline{K}_\e}
\def\Bl{\left(}
\def\Br{\right)}
\def\Blb{\left\{}
\def\Brb{\right\}}
\def\Giv{\ \Big|\ }
\def\tK{{\widetilde K}}
\def\k{\kappa}
\def\re{\R}

\def\tu{{\tilde u}}
\def\tw{{\tilde w}}
\def\ty{{\tilde y}}
\def\tx{{\tilde x}}
\def\tG{{\widetilde G}}
\def\tM{{\widetilde M}}
\def\erm{{\rm e}}
\def\Exp{{\rm Exp\,}}
\def\hx{{\hat x}}
\def\Kln{K_{{\rm ln}}}
\def\Kst{K^*}
\def\jue{_j^\e}
\def\yue{_y^\e}
\def\ue{^\e}
\def\sjN{\sum_{j=1}^N}
\def\non{\nonumber}
\def\adb{}
\def\ccc{c}
\def\bda{b'}
\def\added#1{{\textcolor{blue}{#1}}}
\def\added{}
\begin{document}

\title{General random walk in a random 
    environment defined on Galton--Watson trees}
\author{
A. D. Barbour\footnote{Institut f\"ur Mathematik, Universit\"at Z\"urich,
Winterthurertrasse 190, CH-8057 Z\"URICH;
ADB was supported in part by Australian Research Council Grants Nos DP120102728 and DP120102398.}%\msk}
\ ~and
Andrea Collevecchio\footnote{School of Mathematical Sciences, Monash University, Clayton, VIC 3800, Australia; AC was supported in part by Australian Research Council Grants Nos DP140100559 and ERC  Strep \lq MATHEMACS\rq.}
\\ 
Universit\"at Z\"urich and Monash University
}
\date{}

\maketitle

%\vspace*{0.9cm} 
\noindent {\bf Abstract.}   We consider a particle performing a random walk on a Galton--Watson tree, 
when the probabilities of jumping from a vertex to any one of its neighbours is
determined by a random process. We introduce a method for deriving conditions under
which the walk is either transient or recurrent.
We first suppose that the weights are \iid, and re-prove a result 
of Lyons \& Pemantle~\cite{LyP92}. We then assume a Markovian environment along each line of descent, 
and finally consider a random walk in a Markovian environment that itself changes the environment.  
Our approach involves studying the typical behaviour of the walk on fixed lines of
descent, which we then show determines the behaviour of the process on the whole tree. 

% Our model is a strong generalization of the once-reinforced walk, and it exhibits a phase transition.

\bigskip

%\noindent {\bf AMS 2010 subject classification:} 
\bigskip 


\setcounter{footnote}{0}
\def\thefootnote{\arabic{footnote}}

%\newpage
%\addtocounter{page}{-1}

\section{Introduction}\label{intro}
We consider the behaviour of a random walk in a random environment, which consists of a
randomly sampled Galton--Watson tree, with the jump probabilities at each vertex
% probabilities of jumping from a vertex to any one of its neighbours are
being prescribed by a further random mechanism.  We derive conditions on the environment
under which the walk is transient --- that is, the event that the walk never returns
to the root has positive probability --- and under which it is recurrent, when the
probability of returning to the root is~$1$.
Our approach, which has its roots in that of \cite{Co06}, involves studying the typical 
behaviour of the process on fixed lines of descent, which we then show determines the 
behaviour of the process on the whole tree.  We combine these ideas with suitable 
large deviation principles, and an analysis of the resulting variational formula,
enabling rather satisfactory results to be obtained under relatively weak conditions. 

The  Galton--Watson tree is sampled first,
starting from a root~$\r$. Given the tree, positive weights are then assigned 
at random to its edges, and the jump probabilities are determined from the weights.
In Section~\ref{iid}, we suppose that the weights are \iid, and recover a result 
of Lyons \& Pemantle~\cite{LyP92}, in Theorem~\ref{theorem1}.  
A different proof of the Lyons--Pemantle theorem is contained in \cite{MeP}; see 
also \cite{Far11} for the multitype Galton-Watson case. In Section~\ref{Markov}, 
we extend the argument to an environment in which the values of the weights
evolve as a Markov chain along rays, 
giving sufficient conditions for both transience and recurrence in Theorem~\ref{theorem2}.
Finally, in Section~\ref{one-change}, we illustrate the power of our method
by considering a random walk in a Markovian environment that itself changes the environment;
see Theorems \ref{rmainthrei} and~\ref{rmainthrei2}.  
The results of this section should be 
compared to the behaviour of once-reinforced random walk; see, for example, 
\cite{Co06}, \cite{Dai05} or \cite{Du02}, and the general survey of reinforcement in \cite{Pe07}. 
The model that we discuss is a strong generalization of the once-reinforced walk, and it exhibits multiple 
phases (see Theorem~~\ref{rmainthrei2}), whereas the once-reinforced walk on the supercritical 
Galton--Watson tree is always transient (see \cite{Co06} or \cite{Dai05}).   
 \added{This  should be compared to results obtained in a recent preprint by Kious \&  Sidoravicius  \cite{KS16}.}

Let $\Gcal$ be an infinite tree with root~$\r$.  We augment~$\Gcal$ by adjoining a parent 
$\parent{\r}$ to the root $\r$. If two vertices $\nu$  and $\mu$ are 
the endpoints of the same edge, they are said to be neighbours, and this property is 
denoted by $\nu \sim \mu$. 
The distance $|\nu - \mu|$ between any pair of  vertices $\nu, \mu$, not necessarily adjacent, 
is the number of edges in the unique self-avoiding path connecting $\nu$ to~$\mu$. We set $|\parent{\r}| = -1$. 
For any  other vertex $\nu$,  we let~$|\nu|$ 
be the distance of~$\nu$ from the root~$\r$. 
We denote by~$b(\nu)$ the number of neighbors of~$\nu$ at level $|\nu|+1$, its 
offspring number, and we use~$\parent{\nu}$ to denote the parent of~$\nu$.  We write $\nu < \mu$
if~$\nu$ is an ancestor of~$\mu$.

For $\nu$ a vertex of~$\Gcal$, we write
$$
    \mathbf{A}_{\nu} \Eq (A_{\nu1}, A_{\nu2}, \ldots)
$$
to denote the (finite, positive) weights on the edges between $\nu$ and its offspring. 
\added{For simplicity, we index the weight associated to edge~$e$ by the endpoint of~$e$  with larger distance from~$\r$.}
The environment~$\omega$ for the random walk on the tree
is then defined, for any vertex $\nu$ with offspring $\nu i$, $1\le i \le b(\nu)$,  
by the probabilities
\begin{equation}
  \omega(\nu, \nu i) \Def \frac{A_{\nu i}}{1 + \sum_{1\le j \le b(\nu)} A_{\nu j}}; 
  \qquad \omega(\nu, \parent{\nu} ) \Def \frac 1{1 + \sum_{1\le j \le b(\nu)} A_{\nu j}}.
\end{equation}
We set $\omega(\nu,\mu) = 0 $ if $\mu$ and~$\nu$ are not neighbours.
Given the environment $\omega$, we define the random walk 
$\mathbf{X}= \{ X_{n},\, n \ge 0\}$ that starts at~$\r$ to be the Markov chain with
$\mathbf{P}_{\omega} (X_{0}=\r) =1$, having transition probabilities
$$ 
 \adb{\mathbf{P}_{\omega}(X_{n+1} = \mu_1\;|\; X_{n} = \mu_0) \Eq \omega(\mu_0, \mu_1).}
$$
Moreover, we assume that $\parent{\r}$ is an absorbing state for the walk.
The environment is random in two respects.  First, the Galton--Watson tree~$\Gcal$ is
realized;  then, for each vertex~$\nu \in \Gcal$, the weights~$\mathbf{A}_\nu$ are realized.
The combined probability measure from which the environment is realized is denoted
by~$\mathbb{P}$ and its expectation by $\bbE$, and the semi-direct product 
$\mathbf{P} := \mathbb{P} \times \mathbf{P}_{\omega}$ represents the annealed measure.
The details of the probability measures used to construct the environment are given
in the subsequent sections.

We use  $[\nu, +\infty)$ to denote a generic infinite line of descent from~$\nu$.

\section{Random walks in i.i.d. environment.}\label{iid}
In this section, we assume that~$\Gcal$ is a Galton--Watson tree with offspring mean~$b > 1$.
Given the realization of the tree, we assume that the sets of weights 
$(\bA_\nu,\,\nu\in\Gcal)$ are independent, and that, for each~$\nu$, the weights
$(A_{\nu i},\,1\le i\le b(\nu))$ are exchangeable, with the distributions
of the~$A_{\nu1}$, $\nu\in\Gcal$, all identical.  Under these assumptions, we 
prove the following theorem, first given by Lyons \& Pemantle~\cite{LyP92}, as part of a sharp result.

\begin{theorem}\label{theorem1} 
Assume that  $\Gcal$ and the environment are distributed as above. If 
$\inf_{\lambda \in [0,1]} \mathbb{E}[A_{\r1}^{\lambda}] > b^{-1}$, then $ \mathbf{X}$ is transient; that is, 
with positive probability, $\bX$ does not hit $\parent{\r}$.
\end{theorem}

\nin Our proof relies on the Mogulskii large deviations principle.

We assume that $\bX$ is recurrent and find a  contradiction.
We  consider the behaviour of the random walk $\bX$ observed along any infinite line of descent  
$\s = [\parent{\r}, \infty)$, if one exists. Such lines exist with positive probability, 
since~$b>1$. We call this 
restricted process $\bX\up$. Note that, by our assumption of recurrence, the process $\bX\up$ has the 
following transition probabilities:
\[
    \bP_\o[X_{n+1}\up = \s_{r+1} \giv X_n\up = \s_r] \Eq \frac{A_{\s_{r+1}}}{1 + A_{\s_{r+1}}};\quad
    \bP_\o[X_{n+1}\up = \s_{r-1} \giv X_n\up = \s_r] \Eq \frac{1}{1 + A_{\s_{r+1}}},
\] 
where we denote the successive vertices in~$\s$ by $\s_j$, $j\ge-1$, with $\s_0 := \r$ and $\s_{-1}
:=\parent{\r}$.  We define $T_{-1}$ to be the first 
time~$\mathbf{X}\up$ hits $\parent{\r}$, and $T_{n}$ the first time the process hits~$\s_{n}$. 
Note that the $\bP$-distributions of $T_{-1}$ and $T_{n}$ are not affected by the choice of~$\s$.

\begin{proposition}\label{prop1}
If
\eq\label{ADB-basic-condition} 
   \limsup_{n \ti} \frac 1n  \ln \mathbf{P}(T_{-1} > T_{n}) \ >\ -\ln b,
\en 
then $\mathbf{X}$ is transient.
\end{proposition}

\begin{proofsect}{Proof}  
We mimic the proof in \cite{Co06}. Assume that $\mathbf{X}$ is recurrent. By assumption, 
there exists an~$n^{*}$ such that 
$ b^{n^{*}}\mathbf{P}(T_{-1} > T_{n^{*}}) >1$.  We now construct a branching process as 
follows. Set  $\tau := \inf \{ i > 0 \colon X_{i} = \parent\r\}$.
\ignore{ Here, 
for $\mu\in\Gcal$, $\bX(\mu,m)$ denotes the random walk~$\bX$ restricted to 
$$
  \Gcal(\mu,m) \Def \{\parent\mu,\mu\} \cup \{\nu \in\Gcal\colon\, \nu > \mu,\, |\nu| \le |\mu| + m\}, 
$$
and, if~$\bX$ only spends a finite number of steps in~$\Gcal(\mu,m)$, continued according 
to the Markov chain reflected at each state~$\nu$ such that $|\nu| = |\mu| + m$.}
We color green the  vertices $\nu$ at level $ n^{*}$ which are visited before time $\tau$.
Define 
$$
   S_{\nu} = \inf\{ n\ge 0 \colon X_{n} = \nu\}.
$$ 
Under our assumptions,  $S_{\nu}< \infty$ a.s.\ for each~$\nu$.
A vertex~$\nu$ at level~$j n^{*}$, for some integer $j \ge 2$, is colored green, if its 
ancestor~$\mu$ at level $(j-1) n^{*}$ is green, and $(\bX_{j},\,j \ge S_\mu)$  
hits~$\nu$ before it returns to $\parent{\mu}$. 
The green vertices evolve as a Galton--Watson tree, with offspring mean
$b^{n^{*}}\mathbf{P}(T_{-1} > T_{n^{*}}) > 1$.  Hence this random tree is supercritical, and thus
the probability of there being an infinite number of green vertices is positive.  But this contradicts
the assumption that~$\bX$ is recurrent.
% This, in turn, implies transience.
\ep
\end{proofsect}

\bigskip

\begin{proofsect}{Proof of Theorem \ref{theorem1}}
In view of Proposition~\ref{prop1}, it is enough to show that~\Ref{ADB-basic-condition} is
satisfied.
We use a well-known formula for the hitting probability for random walk in random environment 
(see, for example, Sznitman~\cite{SZ02}, Equation~$44$),
$$ 
   \mathbf{P}(T_{n}< T_{-1}) \Eq \E\Bigl[\bigl(\sum_{r=0}^{n} \prod_{j=1}^{r} A^{-1}_{\s_j}\bigr)^{-1}\Bigr].
$$
Denote by $\floor{x}$ the integer part of~$x$. Then it follows directly,
\adb{because
$$
 \begin{aligned}
    \Bigl(n \max_{r \le n}  \prod_{j=1}^{r} A^{-1}_{\s_j}\Bigr)^{-1} 
      &\le \Bigl(\sum_{r=0}^{n} \prod_{j=1}^{r} A^{-1}_{\s_j}\Bigr)^{-1} 
   \le  \Bigl( \max_{r \le n}  \prod_{j=1}^{r} A^{-1}_{\s_j}\Bigr)^{-1},
 \end{aligned}
$$}
that
\begin{equation}\label{vardh1}
\begin{aligned}
  \liminf_{n \ti}\frac 1n \ln \mathbf{P}(T_{n}< T_{-1}) &\Eq 
          \liminf_{n \ti}\frac 1n \ln \E\Big[\Bigl(\sum_{r=0}^{n} \prod_{j=1}^{r} A^{-1}_{\s_j}\Bigr)^{-1}\Big]\\
  &\Eq  \liminf_{n \ti}\frac 1n \ln \E\Big[ \min_{r \le n} \prod_{j=1}^{r} A^{}_{\s_j}\Big]\\
  &\Eq \liminf_{n \ti} \frac 1n \ln \E\Big[ \erm^{\min_{r \le n} \sum_{j=1}^{r} \ln A_{\s_j}}\Big]\\
  &\ \ge\ \liminf_{n \ti} \frac 1n \ln \E\Big[ \erm^{\min_{t \in [0,1]} \sum_{j=1}^{\nti} \ln A_{\s_j}}\Big].
\end{aligned}
\end{equation}
Denote by $D[0,1]$ the space of  functions $f\colon [0,1] \to \R$, which are right-continuous, 
have limits from the left and have $f(0)=0$. Endow this space with the uniform convergence topology. 
We write~$\mathcal{AC}$ for the subspace of $D[0,1]$ consisting of all absolutely continuous functions. 
By the Mogulskii theorem (see \cite{dembo88}, Theorem 5.1.2), the distribution of 
$\{(1/n)\sum_{j=1}^{\nti} \ln A_{\s_j},\, t \in [0,1]\}$ satisfies a large deviation principle in $D[0,1]$.  
The rate function for this large deviation principle is 
$$ 
  I(f) \Def \int_{0}^{1} \sup_{\lambda} \Bigl\{ f'(t) \lambda - \ln \E[A_{\r1}^{\lambda}] \Bigr\}\, \d t,
$$
if $ f \in \mathcal{AC}$, and $I(f) = +\infty$ if $ f \notin \mathcal{AC}$.
\added{Note that~$I(f)$ is lower semicontinuous, but does not necessarily have compact level sets, so that
it is not necessarily a `good' rate function.}

The function $g \colon \mathcal{AC} \to (-\infty, 0]$ defined by 
$g(f) = \min_{t \in [0,1]} f(t)$ is continuous in $\mathcal{AC}$. 
By  \added{the lower bound in Varadhan's lemma   (see \cite{dembo88}, \added{Lemma 4.3.4}),  we 
%it is sufficient  
%to prove that 
%\begin{equation}\label{v2}
%    \limsup_{n \ti} \frac 1n \ln \E\Big[ \erm^{ 2 \min_{t \in [0,1]} \sum_{j=1}^{\nti} \ln A_{\s_j}}\Big]\ <\ \infty.
%\end{equation} 
%As $ \min_{t \in [0,1]} \sum_{j=1}^{\nti} \ln A_{\s_j} \le 0$,  \eqref{v2} is immediate, and we can apply \added{the lower bound} in
%Varadhan's lemma to
 get}
\begin{equation}\label{lbtran-1}
   \added{\liminf_{n \ti}}\, \frac 1n \ln \E\Bigl[ \erm^{\min_{t \in [0,1]} \sum_{j=1}^{\nti} \ln A_{\s_j}}\Bigr]\ \added{\ge}
    \ \sup_{f \in \mathcal{AC}} \Bigl\{ \min_{t \in [0,1]} f(t) -  I(f) \Bigr\}.
\end{equation}
Since the function $ \phi(\lambda) := \ln \mathbb{E}[A_{\r 1}^{\lambda}]$ is convex, it
follows from Proposition~\ref{varfor1} in the Appendix that the solution to
the variational formula on the right hand side of \eqref{lbtran-1} is given by 
\begin{equation}\label{varformd}
    \sup_{f \in \mathcal{AC}} \Big\{ \min_{t \in [0,1]} f(t) 
      -  \int_{0}^{1} \sup_{\lambda} \{f'(u) \lambda - \ln \E[A_{\r1}^{\lambda}]\}\,  \d u \Big\} 
    \Eq \inf_{\lambda \in [0,1]} \ln \,  \E[A_{\r1}^{\lambda}].
\end{equation}
Combining \Ref{vardh1}, \Ref{lbtran-1} and~\Ref{varformd}, it follows that~\Ref{ADB-basic-condition}
is satisfied, proving Theorem~\ref{theorem1}.
\ep
\end{proofsect}

\ignore{
If the realization of the tree~$\Gcal$ is finite, then the random walk clearly returns to
the root infinitely often with probability one.  If~$\Gcal$ is infinite, then it can be
realized as a {\it spine\/}, consisting of a single line of descent consisting of
individuals having the size-biased offspring distribution, together with independent
copies of~$\Gcal$ emanating from each of the remaining offspring of the individuals
in the spine: see ? Lyons, Pemantle \& Peres~(1995) ?.   Realizing~$\Gcal$ in this
way, and then realizing the weights afterwards, Theorem~\ref{theorem1} shows that, when 
the walk reaches an individual on the spine for the first time, there
is probability at least~$\de$, for some $\de > 0$, that the random walk will next move
to an individual off the spine, and never return to the spine.  Hence, in this case,
the random walk with probability one visits the root only finitely often. 
}

\section{Markovian environment}\label{Markov}

We now show that the proof used in the previous section allows us to treat more general
dependence between the weights, provided that we have a suitable large deviation principle.

Let $\sigma$ be an infinite line of descent $[\r,\infty)$. 
In this section, 
% we analyze the case where the $A_{\s_{i}} \in (0, \infty)$, with $ i \ge 0$, 
we assume that there is 
% a process which satisfies the following. There exists 
a process $\{ M_{\s_{i}}, i \ge 1\}$ in a Polish space~$\Sigma$, such that the pair 
$\Gamma_{\s_i} \Def (A_{\s_{i}}, M_{\s_{i}})$, with $i \ge 0$, is a Markov chain on  
$\S' = (0, \infty) \times \S $, with transition kernel 
\begin{equation}\label{pmc}
  \begin{aligned}
      K(x, B) \Def \mathbb{P}\big(\Gamma_{\s_i} \in   B   \,|\, \Fcal_{i-1} \cap \{\Gamma_{\s_{i-1}} =x\}\big),
  \end{aligned}
\end{equation}
for any  $ B \in \Bcal \Def \Bcal(\S')$; here, $\Fcal_{i}$, $ i \ge 1$, is the natural filtration 
of the process  $\Gamma_{\s_i}$, $i \ge 1$. 

\added{
\begin{remark}   
We assume, for each~$\nu$, that  the random  
variables $(\Gamma_{\nu i},\, i \ge 1)$ are generated  from some joint distribution whose marginals, 
conditionally on $\Gamma_{\nu}$ and the further past, are equal to   $K(\Gamma_{\nu}, \cdot)$. 
Note that we do not need to assume independence among the  $(\Gamma_{\nu i},\, i \ge 1)$. 
The construction can proceed sequentially along the tree, using any initial condition for
$\Gamma_{\r} \in \S'$.
\end{remark} 
}

For any vertex $\nu$, recall that  the set of vertices which are descendants of~$\nu$ consists of 
those vertices~$\mu$ such that~$\nu$ lies on the shortest path connecting~$\mu$ to the root~$\r$. 
We deem~$\nu$ to be its own descendant.  
\ignore{We assume that, given $\Gamma_{\parent\nu}$, the random
element~$\Gamma_{\nu}$ is conditionally independent of the sigma algebra
$$ 
   \sigma\big\{ \Gamma_{\mu} \colon \mu \mbox{ is not a descendant of } \nu \big\}.
$$}
We are motivated by examples where the process $\{A_{\s_i}, i \ge 0\}$ is determined as a functional 
of Markov processes defined on rays.

\ignore{ 
The following very simple example makes it clear that we need some assumption on~$K$ in order to have 
a uniform behaviour, irrespective of the starting point.
\begin{example}  
Suppose that $A_{\s_{i}}$ is a Markov chain on $(0, \infty)$ defined as follows. 
Let $\Dcal_{+}$ be the set of positive 
dyadic numbers, of the form $m/2^{n}$ for some $m,n\in\nat$.  For $x \in \Dcal_{+}$, define 
$K(x, \{x/2\}) =1$. If $x \in (0, \infty)\setminus \Dcal_{+}$, set 
$$ 
     K(x, B) \Def \frac{|B \cap (0, C)|}{C},
$$
where $|\cdot|$ denotes Lebesgue measure and $C$ is a positive constant. If $ A_{\r} \in \Dcal_{+}$, 
then the walk~$\bX$  on the binary tree is recurrent, while, if~$C$ is large enough and $ A_{\r} \notin \Dcal_{+}$, 
then the walk  is transient, by Theorem~\ref{theorem1}.
\end{example}
}

In order to make use of a uniform large deviation principle for Markov
chains, we make the following assumption.  It is somewhat reminiscient of the requirement
for Harris recurrence, but is much stronger, in that many specific measures must be
dominated.  We also make use of the assumption to construct regeneration events
for the environment.

\noindent {\bf Assumption 1.}\ 
{\sl There exist integers $ 0 < \ell \le  N$ and a constant $ \kappa \ge 1$ such that, for all 
$x,y  \in \S'$ and $B \in \Bcal$, we have 
\begin{equation}\label{unifor}
K^{\ssup \ell} (x, B) \le \frac \kappa N \sum_{m=1}^{N} K^{\ssup m}(y, B),
\end{equation} 
where $ K^{\ssup \ell}$ stands for the $\ell$-th convolution of the kernel $K$.
}

%\noindent {\bf Assumption $\mathbf{K(2)}$}\ 
%{\sl We assume that   
%\begin{equation}\label{toone}
%\lim_{\eps \to 0 } \inf_{y \in \S} \mathbb{P}\big(A_{\s_{1}} > \eps \giv  \Gamma_{\s_0} = y\big) = 1.
%\end{equation}
%}

%\begin{remark}
Note that i.i.d.~$\{A_{\s_{i}}\}$ satisfy Assumption~$\ki$, and so does any finite 
state space irreducible Markov chain  
$(A_{\s_{i}}, M_{\s_{i}})$, but there
are of course many other possibilities.

Although the classical results on large deviations require the finiteness of all moments 
(see condition $({\bf\hat{U}})$, page 95 of  \cite{ds89}), 
we do not assume that the support of the~$A_{\s_i}$ is either compact 
or bounded away from zero; nor do we make any assumptions on the moments of~$A_{\s_i}$.
Instead, we use truncation in order to apply the general results.
%\end{remark}
We nonetheless need one further assumption. 
Setting
\begin{equation}\label{defng}
   \eta_{\eps,r} \Def 1 - \inf_{y \in \S'} \mathbb{P}\big(\eps < A_{\s_{1}} \le r \giv  \Gamma_{\s_0} = y\big),
\end{equation} 
we require:

\noindent {\bf Assumption~2.}\ 
{\sl
For $\eta \Def \liminf_{\eps \downarrow 0, r\to\infty} \eta_{\eps,r}$, we have $\eta < 1$.
}

\nin The following example shows that, even when~$A_{\s_i}$ itself is a Markov chain, 
Assumption~$\ki$ does not in general imply Assumption~$\kt$.

\begin{example}  
Suppose that $K(x,\cdot)$ is the mixture $(1-\a)\Exp(1) + \a\Exp(\hx)$, where $\hx := x\vee1$,
$0 \le \a \le 1$ and $\Exp(\l)$ denotes the exponential distribution with mean~$\l^{-1}$.  Then it is easy to
check that $\eta = \a$, and that $K^{(2)}(x,\cdot)$ has a density~$k^{(2)}(x,\cdot)$ satisfying
\[
     (1-\a e^{-1})e^{-w} \Le k^{(2)}(x,w) \Le 3e^{-w},
\]
uniformly in~$x$, so that Assumption~$\ki$ is satisfied with $\ell=2$, but Assumption~$\kt$ is not
satisfied if $\a=1$.
\end{example}

\ignore{
In the example above, the distributions~$K(x,\cdot)$ are continuous in the topology of weak
convergence as~$x$ varies, but mass~$\a$ can be set arbitrarily close to zero by choosing~$x$
large enough.  If this possibility is excluded, and if $\lim_{x\to0}K(x,(0,\e))$ exists
for each $\e > 0$, then Assumption~$\ki$ does imply Assumption~$\kt$ for continuous families
of kernels~$K$, when~$A_{\s_i}$ itself is Markovian.
If some~$M$ is needed to make $(A_{\s_i},M_{\s_i})$ Markov, 
continuity of the family of kernels is not enough for Assumption~$\ki$ to imply Assumption~$\kt$. 
}  
 
For all $x \in \S'$ and for all $ B \in \Bcal(\S')$, define $\S'_\e := (\eps, \infty)\times \S$ and
$$ 
    \oke(x, B)  \Def \frac {K(x, B\cap\S'_\e)}{K(x, (\eps, \infty)\times\S)};
$$
note that~$\oke$ is a probability kernel on~$\S'_\e$, and that it satisfies Assumption~$\ki$
for all $\eps$ such that \adb{$\eta_{\eps,\infty} := 
      1 - \inf_{y \in \S'} \mathbb{P}\bigl(A_{\s_{1}} > \eps \giv \Gamma_{\s_0} = y \bigr) < 1$. }
To prove the latter fact, observe that, for all Borel sets $ B \in \Bcal(\S'_\e)$, 
we have 
\begin{equation}\label{newk0}
  \begin{aligned}
    \oke^{\ssup \ell}(x, B) \Le (1- \eta_{\eps,\infty})^{-\ell} K^{\ssup \ell} (x, B) 
          &\Le \frac \kappa{(1- \eta_{\eps,\infty})^{\ell}N} \sum_{j=1}^{N} K^{\ssup j}(y, B)\\
     & \Le \frac \kappa{(1- \eta_{\eps,\infty})^{\ell}N} \sum_{j=1}^{N} \oke^{\ssup j}(y, B).
  \end{aligned}
\end{equation}

 For any $0 < \e < 1$, and for some $x^{*} \in [1, \infty)\times\S$, define the measure $\mmu_\eps$ 
on~$\S'$ by
\eq\label{ADB-beta-eps-def}
   \mmu_\eps(\cdot) \Def \oke^{\ssup \ell}(x^{*}, \cdot),
\en
where $\ell$ is the same as in Assumption~$\ki$. Set 
$\mmu(\cdot) = \lim_{\e\to0}\mmu_\e(\cdot) = K^{\ssup\ell}(x^*,\cdot)$.

\begin{proposition}\label{prop2}
Under Assumptions $\ki$ and~$\kt$, if
\begin{equation}\label{condmar}
 \liminf_{\e \to 0}  \liminf_{n \ti}  \frac1n \int_{\S'} \ln  \bP\bigl(T_{-1} > T_{n} \;|\; \Gamma_{\r}= y \bigr) 
       \mmu_{\e}(\d y) 
             \ >\  -\ln b,
\end{equation}
then $\mathbf{X}$ is transient.
\end{proposition}

\begin{proofsect}{Proof}
%Fix $c>0$, to be chosen later. In view of \eqref{condmar}, and using Jensen's inequality, we can 
%find~$n^{*}= n^{*}(c) $ such that 
%% for any $m \ge n^{*}$ we have 
%\begin{equation}\label{condmar0}
%   b^{n^*} \int_{\S'} \bP(T_{-1} > T_{n^*}  \,|\, \Gamma_{\r} = y) \mmu(\d y) \ >\ 1/c.
%\end{equation}
%Since $\mmu = \lim_{\e\to0}\mmu_\e$, this implies that we can choose $\eps >0$ small enough that
%\begin{equation}\label{condmar2}
%   b^{n^*} \int_{\S'} \bP(T_{-1} > T_{n^*}  \,|\, \Gamma_{\r} = y) \mmu_{\eps}(\d y) \ >\ 1/c.
%\end{equation}
%\ignore{
%where, for $y = (u,z)$, 
%\begin{equation}\label{mmab}
%   \frac{\d \mmu_\eps^*}{\d\mmu_\eps}(y) \Eq c_{\mmu,\e}\,\frac{u}{1+u}.
% \end{equation}  
%}
Under the assumption that~$\mathbf{X}$ is recurrent,
we construct a random subtree of $\Gcal$, consisting  of green 
vertices, that contains a number of vertices stochastically larger than the number of vertices in   
a supercritical Galton--Watson tree. These green vertices are such that the random
walk~$\bX$ visits them before it first reaches~$\parent{\r}$.
The fact that this random subtree is infinite with positive 
probability implies a contradiction, and hence that~$\mathbf{X}$ is transient. 

A direct calculation shows that, for any $y\in\S$ and $1\le j\le N$,
\eqs
    p\jue(y,B) &:=& \bbP\Bl \bigcap_{l=1}^{j-1} \{A_l \ge \e\},\,\G_j \in B \Giv \G_0=y \Br \\
              &\ge& (1-\eta_{\eps,\infty})^N \oke^{\ssup j}(y,B),
\ens
for all $ B \in \Bcal(\S'_\e)$.
It thus follows, from \Ref{newk0} and~\Ref{ADB-beta-eps-def}, that
if~$U$ is uniformly distributed on $\{1,2,\ldots,N\}$, 
independently of~$\G$, then, for all $y\in\S'$ and $ B \in \Bcal(\S'_\e)$,
\eqa
   p\yue(B) &:=& \frac1N \sjN p\jue(y,B) 
   \Eq \bbP\Bl \bigcap_{l=1}^{U-1} \{A_l \ge \e\},\,\G_U \in B \Giv \G_0=y \Br \non\\
            &\ge& \k^{-1}(1-\eta_{\eps,\infty})^{N+\ell}\oke^{\ssup\ell}(x^*,B) \ =:\ \de_\e \mmu_\e(B), 
    \label{ADB-lower-measure}
\ena
with $\de_\e > 0$ for all~$\eps$ small enough, since $\eta_{\eps,\infty} < 1$ for all~$\eps$ small
enough, in view of Assumption~$\kt$.
Because of~\Ref{ADB-lower-measure}, $\de_\e \mmu_\e$ is absolutely continuous with respect to~$p\yue$,
and
\[
   \de_\e \,\frac{d\mmu_\e}{dp\yue}(y') \Le 1 \Eq \sjN f\jue(y,y'),
\]
where
\[
    f\jue(y,y') \Def \frac1N \frac{dp\jue(y,\cdot)}{dp\yue}(y'),\quad 1\le j\le N.
\]
Hence, if we set  $g\jue(y,\cdot) := \de_\e \frac{d\mmu_\e}{dp\yue}(y') f\jue(y,y')$,  $1\le j\le N$,  
it follows that $0 \le g\jue(y,y') \le f\jue(y,y')$ 
for all~$y' \in \S'_\e$ and  $\sjN g\jue(y,y') =  \de_\e \frac{d\mmu_\e}{dp\yue}(y')$.
 
This justifies the following construction.  Starting at a vertex~$\nu$ that has an infinite
line of descent, let~$\G_0$ denote the value~$y\in\S'_\e$ at~$\nu$. Realize~$U = U_\nu$ 
uniformly distributed on $\{1,2,\ldots,N\}$ and a random variable~$U'$ uniformly distributed on~$[0,1]$,
independently of all else.  Because there is an infinite line of descent from~$\nu$, there is at
least one line of descent from~$\nu$ of length~$U$; if there is more than one, choose one
at random.  Denote it by $\nu_1,\ldots,\nu_U$, and set $\nu_0 := \nu$.
Independently, realize the chain~$\G$ along this line of descent, starting from~$\G_0$ at~$\nu$.
Say that the event~$E_\nu$ occurs if $A_j \ge \eps$, $1\le j\le U-1$, and if 
$U' f_U\ue(y,\G_U) \le g_U\ue(y,\G_U)$. In this way,
the distribution~$\mmu_\e$ is obtained as the distribution of~$\G_U$ on an event~$E_\nu$ of
probability~$\de_\e$, and with $A_j \ge \e$, $1\le j\le U$.
For any pair of vertices $\nu$, $\mu$, with $\mu$ a descendant of $\nu$,  denote by 
$\bX(\nu, \mu)$ the process $\bX$ restricted to the finite graph consisting of the vertices 
in the finite ray $[\parent{\nu}, \mu]$ and the edges connecting them.
A vertex~$\nu'$ is green if it has an infinite line of descent, and 
is descended from a green vertex~$\nu$ in the following way.
% For $\nu$ to be a green vertex, 
% It is good 
$E_\nu$ must occur, and \adbc{then~$\bX(\nu,\nu_U)$} has to reach~$\nu_U$ before hitting~$\nu^{-1}$;
the latter event has probability at least~$\{\e/(1+\e)\}^N$.  Finally, $\nu'$ should be a
descendant at distance~$n$ from~$\nu_U$,
and \adbc{$\bX(\nu_U,\nu')$} should reach~$\nu'$  before it 
hits~$\parent\nu_U$, an event of probability
\[
   \int_{\S'_\e} \bP( T_{-1} > T_{n} \giv \G_\r = y)\mmu_\e(dy).
\]
% and if it has \added{an} infinite line of descent, an event of conditional probability \hbox{$1-q$.}
Thus the expected number of green `offspring' of a green vertex~$\nu$ is at least
\begin{equation}\label{alln}
     b^n \de_{\e} \{\e/(1+\e)\}^N \Blb \int_{\S'_\e} \bP( T_{-1} > T_{n} \giv \G_\r = y)\mmu_{\e}(dy)
       \Brb (1-q),
\end{equation}
\adbc{where~$q$}  denotes the probability of the extinction of the underlying Galton--Watson tree.
Next, we show that  \eqref{condmar} implies that we can choose $\e$ small enough and $n$ large enough 
that the quantity in \eqref{alln} becomes larger than 1. By taking the natural logarithm of~\eqref{alln} 
and dividing by $n$, we have 
\begin{equation}\label{cm2}
    \frac 1n \ln  \left(\de_{\e} \{\e/(1+\e)\}^N (1-q)\right) + \ln b 
           + \frac 1n\ln  \int_{\S'_\e} \bP( T_{-1} > T_{n} \giv \G_\r = y)\mmu_{\e}(dy).
\end{equation}
Fix $\e>0$ such that 
$$ 
   \liminf_{n \ti}  \frac1n \int_{\S'_\e} \ln  \bP\bigl(T_{-1} > T_{n} \;|\; \Gamma_{\r}= y \bigr) 
           \mmu_{\e}(\d y) \ >\ -\ln b,
$$
as we may, in view of~\eqref{condmar}.  Then, for this choice of~$\e$,
the liminf of~\eqref{cm2} as $n \ti$ is larger than~$0$, proving that the quantity 
in~\eqref{alln} is larger than~$1$  for our choice of $\e$ and for~$n$ large enough.

   By construction, the
distribution of the number of green offspring is the same for all green vertices.  Hence, choosing 
an appropriate~$\e$, and then~$n$ large enough that the quantity in~\eqref{alln} is larger than one,
the Galton--Watson tree of green vertices is supercritical.
\ep
\end{proofsect}

\medskip
The proofs that follow rely on large deviations results. These cannot be directly applied to~$A$, 
so we need to consider truncations. 
For this reason, it is  convenient to introduce the large deviations results that we shall
use applied to a generic process~$W := (W_i,\,i\ge0)$,  which, together with a process $\widetilde{M}$ on~$\S$, 
makes~$\widetilde{\Gamma}$ defined by $\widetilde{\Gamma}_i \Def (W_i, \widetilde{M}_i)$ 
 a Markov chain on~$\re\times\S$.  Let~$\tK$ denote the kernel of this process.

Define  
\begin{equation}\label{lambdadef}
  \L^{\ssup \tK}(\lambda) \Def \limsup_{n \ti}\sup_{\ty\in\re\times\S} \frac 1n \ln \mathbb{E} \big[  
            \erm^{\lambda \sum_{i=1}^{n} W_{i}} \giv  \widetilde{\Gamma}_{\r} = \ty \big], 
   \quad\ \L^{*}_{\tK} (x) \Def \sup_{\lambda}\{ \lambda x - \L^{\ssup \tK}(\lambda)\},
\end{equation}
and let 
$$ 
  S^{\ssup {\tK}}_{n}(t) \Def \frac1n \sum_{j=1}^{\nti}  W_{i}, \qquad  t \in [0,1].
$$

\begin{theorem}\label{LDPL} Fix $0 < C< R < \infty$, and assume that $W_i \in (C, R)$ a.s., for each $i$.    
If $\tK$ satisfies Assumption~$\ki$, then, for any $\Theta \in \Bcal^{+}$, we have 
$$
\begin{aligned}
 -\inf_{x \in \Theta^{\circ}}\Lambda^{*}_{\tK}(x) &\Le \liminf_{n \ti} \frac 1n 
          \ln \inf_{y \in \S'}\mathbb{P}(S_{n}^{\ssup \tK}(1) \in \Theta \;|\; \widetilde{\Gamma}_{\r} = y) \\ 
    &\Le \limsup_{n \ti} \frac 1n \ln \sup_{y \in \S'}\mathbb{P}(S_{n}^{\ssup \tK}(1) \in \Theta 
          \;|\;  \widetilde{\Gamma}_{\r}= y) \Le - \inf_{x \in \overline{\Theta}} \Lambda_{\tK}^{*}(x).
\end{aligned}          
$$
\end{theorem}

\ignore{
\begin{remark}
There are examples where the kernel of the  process $\Gamma_{\s_{i}}$ satisfies $\ki$ and $ \eta >0$. 
In fact, fix $\alpha \in (0,1)$ and suppose that  $ K(2^{-n}, \cdot)$, for $n \ge 1$, has density equal 
to $\alpha 2^{n+1}$ in the interval $I_{n} = (2^{-(n+1)} - 2^{-(n+2)}, 2^{-(n+1)} + 2^{-(n+2)})$, $1-\alpha$ 
in the interval $ [0,1]\setminus I_{n}$, and $0$ elsewhere. For $ x \notin \{2^{-n}, n \ge 1\}$, let 
$ K(x, \cdot)$ have density $1$ in the interval $[0, 1]$ and $0$ otherwise.
Notice that for any borel subset $B$ of $[0, 1]$, and $j \ge 2$, we have that  $ K^{\ssup j}(y, B) = |B|$, 
for all $y \in [0,1]$. Hence,  for any $y \in [0,1]$, 
$$ \sum_{j=1}^{3} K^{\ssup j}(y, B) \ge 2 |B|.$$
This implies that
$$ K^{\ssup 2}(x, B) \le  \frac{1}2 \sum_{j=1}^{3} K^{\ssup j}(y, B).$$
 Hence for the choice of  $ \ell =2$,  $M = 3/2$ and $N=3$ we have that  $K$ satisfies $\ki$. 
On the other hand, for any $\eps >0$, for all $n$ large enough
 $$  1 -  \bbP(A_{\s_{1}} > \eps \giv A_{\s} = 2^{-n}) \ge \alpha.$$
 Hence 
 $$ \eta_{\eps} = \sup_{y \in [0,1]} 1 -  \bbP(A_{\s_{1}} > \eps \giv A_{\s} = y) \ge \alpha,$$
 which implies that $\eta \ge \alpha$.
\end{remark}
}

\begin{proofsect}{Proof of Theorem \ref{LDPL}}  The kernel $\tK $ satisfies condition $({\bf\hat{U}})$, 
page 95 of  \cite{ds89}. Hence,   
the theorem is a consequence of the more general Theorem~4.1.14, page 97 of \cite{ds89},  combined with (4.1.24) 
page 100 of \cite{ds89}, to identify the rate function.
\qed
\end{proofsect}

\medskip
Recall that $ \mathcal{AC}$ denotes the space of absolutely continuous functions~$f$ defined on $[0,1]$, 
with $f(0)=0$ 
and $D[0,1]$ the space of  functions $f$ which are right continuous and have limits from the left,  and have $f(0)=0$.
Both spaces are endowed with the uniform convergence topology. The following result is due to Dembo \& Zajic  \cite{dembo-zajic}.

\begin{theorem}\label{dem}  
Under the hypotheses of Theorem~\ref{LDPL}, the sequence $\{S^{\ssup \tK}_{n}(t), t \in [0,1]\}$ in $D[0,1]$ 
satisfies a large deviations principle with the good, convex, rate function 
$$ 
   I^{*}_{\tK}(f) \Def  
      \begin{cases}
             \int_{0}^{1} \L^{*}_{\tK}(\dot{f}(u)) \,\d u, &\mbox{if $ f \in \mathcal{AC}$}\\
                +\infty, &\mbox{otherwise}.
      \end{cases}
$$
\end{theorem} 

\begin{proofsect}{Proof} 
In virtue of Theorem~\ref{LDPL},  $S^{\ssup \tK}_{n}(1)$ satisfies a uniform large deviations principle. 
We can then use
Dembo \& Zajic (\cite{dembo-zajic}, Theorem~3a) to conclude that $\{S_{n}^{\ssup \tK}(t),\, t\in [0,1]\}$ 
satisfies an LDP with rate function~$I^*_{\tK}(\cdot)$.
\ep	
\end{proofsect}

\medskip
Note that, since $(A_{\s_i}, M_{\s_i})$ is a Markov chain in $(0, \infty)\times\S$, then 
$(\ln A_{\s_i}, M_{\s_i})$ is a Markov chain on~$\R\times\S$.  Define the kernel
$$ 
   \Kln ((\tu,z),B) \Def K((\erm^{\tu},z), E(B)), \qquad \tu \in \R,\ z\in\S,\ B \in \Bcal,
$$
where $E(B) := \{(e^{\tu},z)\colon\,(\tu,z)\in B\}$.  Note that, if $K$ satisfies Assumption~$\ki$, 
then so does the kernel~$\Kln$.

For  $R  \in (0, \infty)$ and $ C \in [-\infty, 0)$,  define the probability kernel $Q_{C, R}$ 
on $(C,R] \times \S$ by
\begin{equation}\label{defnq}
    Q_{C, R}\big(\ty, (\d \tu, \d z)\big) 
             \Def \frac{\Kln \big(\ty, (\d \tu, \d z)\big)}{\Kln(\ty, (C,  R]\times \S)},
\end{equation}
and set $Q_{R} := Q_{-\infty, R}$.

% \qquad  \Lambda^{*}(x) \df   \sup_{\lambda} \lambda x -   \Lambda(\lambda),$$
 
\begin{theorem}\label{theorem2} 
If~$K$ satisfies Assumption~$\ki$, then 
\begin{itemize}
\item[{\rm (i)}] If Assumption~$\kt$ holds, the condition
\begin{equation}\label{condd10}  
    \limsup_{\min\{-C,R\}\to\infty}\,\inf_{\lambda \in [0,1]} \Lambda^{\ssup {Q_{C, R}}}(\lambda) 
             \ >\ -\ln b - \ln (1- \eta) 
\end{equation}
implies transience of $\mathbf{X}$ on $\Gcal$. {The constant $\eta$ is \adbc{the one} introduced 
in Assumption~2.}
\item[{\rm (ii)}] The condition
\begin{equation}\label{con11} 
     \inf_{\lambda \in [0,1]} \Lambda^{\ssup {\Kln}}(\lambda)\ <\ -\ln b 
\end{equation}
 implies positive recurrence of $\mathbf{X}$ on $\Gcal$.
\end{itemize}
\end{theorem}

\begin{remark}
   \adb{In the case of an i.i.d.\ environment, \eqref{condd10}  coincides with the condition 
$\inf_{\lambda \in [0,1]} \Lambda^{\ssup {\Kln}}(\lambda)\ > \ -\ln b$, which is then also
the same as that of Theorem~{\rm \ref{theorem1}}.} 
\end{remark}

\begin{corollary}\label{impco}
 In the uniformly elliptic case, i.e. if  there exists $\eps >0$ such that 
$$ 
 \inf_{x} K\big(x, (\eps, \eps^{-1})\times \Sigma\big) = 1,
$$
\adbc{so that then~$\eta = 0$, we have the following sharp transition:}
\begin{itemize}
\item[{\rm (i)}] The condition
\begin{equation} 
     \inf_{\lambda \in [0,1]} \Lambda^{\ssup {\Kln}}(\lambda)\ >\ -\ln b 
\end{equation}
implies transience of $\mathbf{X}$ on $\Gcal$. 
\item[{\rm (ii)}] The condition
\begin{equation} 
     \inf_{\lambda \in [0,1]} \Lambda^{\ssup {\Kln}}(\lambda)\ <\ -\ln b 
\end{equation}
 implies positive recurrence of $\mathbf{X}$ on $\Gcal$.
\end{itemize}
\end{corollary}

The following examples show particular ways to compute bounds for $\Lambda^{\ssup {\Kln}}$.

\begin{example}\label{mainex}
Suppose that~$A_{\s_i}$  evolves as a discrete irreducible aperiodic Markov chain, 
with state space $ \Xi= (a_1,  a_2, \ldots, a_\ell)$, where $a_{i} \in (0 , \infty)$ for all $i$,  
and with transition matrix $K= (k_{i,j},\, 1\le i,j \le \ell)$.
% Moreover suppose there exists two probability  distributions  $\vec{p}= \{p_i, i \le \ell\}$ 
% and $\vec{q}=\{q_i, i \le \ell\}$ both on $\Xi$ such that each raw of $K$ is stochastically larger 
% that $\vec{q}$ and stochastically smaller than $\vec{p}$. 
 Note that, in the finite case, $\Lambda^{\ssup {\Kln}}(\lambda)$ coincides with  $\ln \rho(\lambda)$, 
where $\rho(\lambda)$ is the  Perron-Frobenius eigenvalue of the matrix whose $(i,j)$th entry 
is $a_{j}^\lambda k_{i,j}$, (see \cite{dembo88}, page~74).
Using the Gershgorin circle theorem, the Perron--Frobenius eigenvalue is bounded above by the 
\adbc{largest}  row sum.
Hence 
$$ 
     \rho(\lambda)\ \le\ \adbc{\max_{i}} \sum_{j=1}^\ell k_{i,j} a_{j}^\lambda.
$$
Hence, Corollary \ref{impco} implies that if  
$$  
     \adbc{\inf_{\l \in [0,1]} \max_{i}} \sum_{j=1}^\ell k_{i,j} a_{j}^\lambda\ <\ 1/b,
$$
then the process is recurrent. The
Gershgorin circle theorem can also be used to get the lower bound 
$$ 
     \rho(\lambda)\ \ge\ \adbc{\min_{i}} \left(k_{i,i} a_i^\lambda  
      - \sum_{j \colon j \neq i}  k_{i, j}  a_j^\lambda\right),
$$
useful if~$K$ is close to being diagonal.
Thus Corollary~\ref{impco} implies that if  
$$
  \adbc{\inf_{\l \in [0,1]}\min_{i}} \left(k_{i,i} a_i^\lambda  
          - \sum_{j \colon j \neq i}  k_{i, j}  a_j^\lambda\right)\ > \ 1/b,
$$
then the process is transient.
\end{example}

In the case $\ell=2$, $\Lambda^{\ssup {\Kln}}(\l)$ can of course be computed
explicitly.  For 
$$
K \ :=\ \left(
   \begin{array}{cc}
      \alpha & 1 - \alpha\\
      1 - \beta & \beta
   \end{array}
  \right),
$$
with $\alpha, \beta \in (0, 1)$, we have
\begin{equation}\label{eigtwo}
   \L^{\ssup {\Kln}}(\lambda) = \ln \frac 12 \Bigl(\alpha a_{1}^{\lambda} + \beta a_{2}^{\lambda}
   + \sqrt{ (\alpha a_{1}^{\lambda} + \beta a_{2}^{\lambda})^{2}
            + 4(1 - \alpha - \beta)a_{1}^{\lambda}a_{2}^{\lambda}  }\Bigr).
\end{equation}

The procedure used in Example~\ref{mainex} can be carried out in continuous space through discretization, 
as the following simple example shows.

\begin{example}\label{cts-ex} 
Let~$K$ be the kernel of a Markov process with compact state space $U$.  
Consider  a finite cover of $U$, say $U_{1}, U_{2}, \ldots, U_{\ell}$, with the property 
that if $i \neq j$ then $U_i \not\subset U_j$, and an $\ell\times \ell$  matrix with 
strictly positive elements $B = \{ b_{i,j}\}$ such that, for all $t>0$, 
$$ 
   K(y, (0, t])\ \le\ \sum_{j=1}^{\floor{t}} b_{i, j} \qquad \forall { y \in U_{i}}.
$$
We emphasize that $B$ need not be a transition matrix.  Set  $\Xi= \{a_{i}, i \le \ell\}$, 
where $a_{i} = \sup \{ b \colon b \in U_{i}\}$. Then
$$ 
   \L^{\ssup {\Kln}}(\lambda)\ \le\ \limsup_{n \ti} \sup_{j_0 \le \ell} 
    \frac 1n \ln \sum_{j_1 =1 }^\ell \ldots \sum_{j_n =1 }^\ell 
         {\rm e}^{\lambda \ln a_{j_1}} \cdots   {\rm e}^{\lambda \ln a_{j_n}} 
           b_{j_0, j_1}\cdots b_{j_{n-1}, j_n}.
$$
Hence, the  logarithm of the Perron-Frobenius eigenvalue of the matrix $\{a_{j}^{\lambda} b_{i,j}\}$ 
is an upper bound for $\L^{\ssup{\Kln}}$.  We can then proceed as in the previous example to determine 
a sufficient condition for recurrence.
An analogous procedure, with lower bounds, can be applied to derive sufficient conditions for transience.
\end{example}

\begin{proofsect}{Proof of Theorem~\ref{theorem2}}
We first prove that condition \eqref{condd10} implies that~\Ref{condmar} holds, and hence, 
by Proposition~\ref{prop2}, that~$\bX$ is transient. Observe that, for any $r>0$, we have
\begin{equation}\label{con2}
\begin{aligned}
  {\liminf_{\e \to 0}}  \liminf_{n \ti} & \int_{\S'}  \frac 1n \ln \bP(T_{n}< T_{-1} \;|\; 
                    \Gamma_{\r} = y) \mmu_{{\e}}(\d y)\\
    &=  {\liminf_{\e \to 0}}   \liminf_{n \ti}  \int_{\S'} \frac 1n 
       \ln \E\Big[\big(\sum_{l=0}^{n} \prod_{j=1}^{l-1} A^{-1}_{\s_j}\big)^{-1} 
                         \;|\; \Gamma_{\r} = y \Big] \mmu_{{\e}}(\d y) \\
    &\ge {\liminf_{\e \to 0}}  \liminf_{n \ti}  \int_{\S'} \frac 1n 
       \ln \E\Big[\big(\sum_{l=0}^{n} \prod_{j=1}^{l-1} (A_{\s_j}\wedge r)^{-1}\big)^{-1} 
                         \;|\; \Gamma_{\r} = y \Big] \mmu_{{\e}}(\d y) \\
    &\ge \liminf_{n \ti} \inf_{y \in \S'} \frac 1n 
       \ln \E\Big[\big(\sum_{l=0}^{n} \prod_{j=1}^{l-1} (A_{\s_j}\wedge r)^{-1}\big)^{-1} 
                         \;|\; \Gamma_{\r} = y \Big]\\
    &\ge\   \liminf_{n \ti} \inf_{y \in \S'}\frac 1n \ln 
                \E\Big[ \erm^{ \min_{t \in [0,1]} \sum_{i=1}^{\nti}\ln (A_{\s_i}\wedge r)} 
            \;|\; \Gamma_{\r} = y \Big].
\end{aligned}
\end{equation}
For $r > 1>\ccc >0$,  let $A_{\s_{i}}(\ccc, r) = (A_{\s_{i}}\vee\ccc) \wedge r$, and set $ C = \ln \ccc$ 
and $ R = \ln r$. Then, writing $\ty_{j} := (\tu_{j},z_{j}) \in \re\times\S$ for $j \ge 1$ and 
$\ty_0 := (\ln u,z)$ for $(u,z) = y$, we have 
\begin{equation}\label{conditimp}
 \begin{aligned}
  &\liminf_{n \ti} \inf_{y \in \S'} \frac 1n \ln 
             \E\Big[ \erm^{\min_{t \in [0,1]} \sum_{j=1}^{\nti} \ln (A_{\s_j}\wedge r)} \giv \Gamma_\r = y\Big]\\
  &\ \ge\ \liminf_{n \ti} \inf_{y \in \S'} \frac 1n \ln 
             \E\Big[ \erm^{\min_{t \in [0,1]} \sum_{j=1}^{\nti} \ln (A_{\s_j}\wedge r)} 
                    \1_{\bigcap_{i=1}^{n} \{A_{\sigma_{i}}> \ccc\} } \giv \Gamma_\r = y\Big]\\
  &\ =\ \liminf_{n \ti} \inf_{y \in \S'} \frac 1n \ln 
             \E\Big[ \erm^{\min_{t \in [0,1]} \sum_{j=1}^{\nti} \ln A_{\s_j}(\ccc, r)} 
                    \1_{\bigcap_{i=1}^{n} \{A_{\sigma_{i}}> \ccc\} } \giv \Gamma_\r = y\Big]\\
  &\ \ge\ \liminf_{n \ti} \inf_{y \in \S'} \int_{([C, R]\times \S)^{n}}   
          \erm^{\min_{t \in [0,1]} \sum_{j=1}^{\nti}  \tu_{j}} 
             \prod_{j=1}^{n} \Kln \big(\ty_{{j-1}}, \d \ty_{j}\big) \,.  
 \end{aligned}
\end{equation}
\adb{Choosing~$\ccc$ small enough and~$r$ large enough that 
\[  
   \inf_{y \in \S'} \Kln(y, [\ccc, r]\times\S)\ \ge\ (1- \eta_{\ccc,r})\ >\ 0,
\]
as we may, because $\eta < 1$,}  we have
\eq\label{new-star}
\begin{aligned}
  \liminf_{n \ti} \inf_{y \in \S'} &\frac 1n \ln \E\Big[ \erm^{\min_{t \in [0,1]} \sum_{j=1}^{\nti} (\ln A_{\s_j}\wedge r)}
          \giv \Gamma_\r = y\Big]\\
  & \ge \liminf_{n \ti} \inf_{y \in \S'} \frac 1n \ln \Big[(1- \eta_{\ccc,r})^{n} \int_{([C, R]\times \S)^{n}}  
        \erm^{\min_{t \in [0,1]} \sum_{j=1}^{\nti}  \tu_j} \prod_{j=1}^{n} 
            \frac{\Kln \big(\ty_{j-1}, \d \ty_j)\big)}{\Kln(\ty_{j-1}, (C, R]\times \S)} \Big]  \\
  &=\ \liminf_{n \ti} \inf_{y \in \S'} \frac 1n \ln \Big[(1- \eta_{\ccc,r})^{n} 
       \widetilde {\mathbb{E}} \big[\erm^{\min_{t \in [0,1]} \sum_{j=1}^{\nti} W_j} \giv \Gamma_\r = y\big] \Big]\\
  &=\ \ln (1- \eta_{\ccc,r}) + \liminf_{n \ti} \inf_{y \in \S'} \frac 1n 
       \ln \widetilde {\mathbb{E}} \big[\erm^{\min_{t \in [0,1]} \sum_{j=1}^{\nti} W_j} \giv \Gamma_\r = y \big]\,,
\end{aligned}
\en
where $\widetilde {\mathbb{E}}$ is the expectation with respect to the Markov chain~$\tG = (W,\tM)$ with probability kernel 
$Q_{C, R}\big(\ty', \d\ty\big)$ introduced in~\eqref{defnq}.

Next we prove that the kernel $Q_{C, R}$  satisfies Assumption~$\ki$.
Note that, for Borel sets  $F \subset (C, R]$ and $E \in \S$, and for any $\tx,\ty \in (C,R]\times\S$, we have  
\begin{equation}\label{newk1}
  \begin{aligned}
  Q_{C, R}^{\ssup \ell} (\tx, F \times E) \Le (1- \eta_{\ccc,r})^{-\ell} \Kln^{\ssup \ell} (\tx, F \times E) 
            &\Le \frac M{(1- \eta_{\ccc,r})^{\ell}N} \sum_{j=1}^{N} \Kln^{\ssup j}(\ty, F \times E)  \\
  & \Le \frac M{(1- \eta_{\ccc,r})^{\ell}N} \sum_{j=1}^{N} Q_{C, R}^{\ssup j}(\ty, F \times E).
  \end{aligned}
\end{equation}
In the last step,  we have used the inequality
$$ 
   \Kln^{\ssup n}(\ty, F \times E) \le Q_{C, R}^{\ssup n}(\ty, F \times E),
$$
valid for $F \subset (C, R]$ and $n \ge 1$, which is easily proved by induction.

Combining Theorem~\ref{dem}  with  Varadhan's lemma, using the {\bf uniform} large deviations 
stated in Theorem~\ref{LDPL},  we find that
\begin{equation}\label{lbtran}
 \begin{aligned}
  \liminf_{n \ti}  \inf_{y \in \S'} \frac 1n 
     &\ln \E\Big[ \erm^{\min_{t \in [0,1]} \sum_{j=1}^{\nti} (\ln A_{\s_j}\wedge r)}
         \giv \Gamma_\r = y\Big]\\
    &\ge\ \ln (1- \eta_{\ccc,r}) +
      \sup_{f \in \mathcal{AC}} \big\{ \min_{t \in [0,1]} f(t) -  I^{*}_{Q_{C, R}}(f) \big\},
 \end{aligned}
\end{equation}
and, since the function $\Lambda^{\ssup {Q_{C, R}}}$ is convex for any $C$ and~$R$, 
Proposition~\ref{varfor1} can be used to solve
the variational formula on the right hand side of~\eqref{lbtran}, giving 
\begin{equation}\label{varformd-2}
\begin{aligned}
    \liminf_{n \ti}  \inf_{y \in \S'} \frac 1n 
            \ln \E\Big[ \erm^{\min_{t \in [0,1]} \sum_{j=1}^{\nti} (\ln A_{\s_j}\wedge r)}
     \giv \Gamma_\r = y\Big] &\ 
       { \ge\  \ln (1- \eta_{\ccc,r}) + \inf_{t \in [0,1]} \Lambda^{\ssup {Q_{\ln\ccc, \ln r}}}(t). }
%    &\ \ge\ \ln (1- \eta_{\eps}) + \inf_{t \in [0,1]} \Lambda^{\ssup {Q_{R}}}(t).
\end{aligned}
\end{equation}
% In the last step, we used the fact that  $ \Lambda^{\ssup {Q_{C, R}}}$ is decreasing in~$C$, and we 
% compared it with $ \Lambda^{\ssup {Q_{-\infty, R}}} = \Lambda^{\ssup {Q_{R}}}$.

Recalling~\Ref{con2}, we thus have
\[
   \liminf_{\e\to0} \liminf_{n \ti} \int_{\S'_\e}  \frac 1n \ln \bP(T_{n}< T_{-1} \;|\; \Gamma_{\r} = y) 
        \mmu_\e(\d y) 
           \ \ge\ \ln (1- \eta_{\ccc,r}) + \inf_{t \in [0,1]} \Lambda^{\ssup {Q_{\ln\ccc, \ln r}}}(t),
\]
{ for any $\ccc,r > 0$ such that $\eta_{\ccc,r} < 1$.
By letting $\ccc\to0$ and~$r\to\infty$,
% taking the supremum over~$R$, which coincides with the limit as $R \uparrow \infty$, 
we get
\eqs
   \lefteqn{\liminf_{\e\to0} \liminf_{n \ti} \int_{\S'_\e}  \frac 1n \ln \bP(T_{n}< T_{-1} \;|\; \Gamma_{\r} = y) 
                \mmu_\e(\d y)}\non\\
   &&\ge\ \ln (1- \eta) + \limsup_{\min\{1/\ccc, r\} \to \infty} \inf_{t \in [0,1]} 
           \Lambda^{\ssup {Q_{\ln \ccc,\ln r}}}(t) \ >\ -\ln b,
\ens}
using~\Ref{condd10}, and~(i) follows from Proposition~\ref{prop2}.
 
Next, we prove that if~\Ref{con11} holds,
then the process is positive recurrent. In this case we just mimic the proof by Lyons \&~Pemantle (see \cite{LyP92}, 
proof  of Theorem~1.3, page 130). We include the proof for sake of completeness and clarity.
 \added{We use the well known  fact {(see \cite{KSK} Proposition 9-131)} that the random walk is positive recurrent if the  
sum of conductances is a.s.\ finite, i.e.
\begin{equation}\label{sucofi} 
    \sum_{\nu\in\Gcal} \prod_{i=1}^{|\nu|} A_{\nu^{-i}}\ <\ \infty \qquad \mbox{a.s.}
\end{equation}
}
From \Ref{con11} and the definition of $\L$, we can choose $t_{0} \in (0,1]$ such that 
$$  
       \mathbb{E}\Bigl[\exp\Bigl\{ t_{0} \sum_{i=1}^{n} \ln A_{i}\Bigr\} \Giv \tG_0 = \ty\Bigr]\ < \ (1/\bda)^n , 
$$
for some $\bda>b$, for all~$\ty$ and all $n$ large enough. \added{Recall that $\nu^{-i}$ denotes the $i$-th ancestor of $\nu$.} Because the branching number of the Galton--Watson tree 
is~$b$, this implies that
\eq\label{ADB-LP1}
     \bbE\Bl \sum_{\nu\colon |\nu|=n} \prod_{i=1}^n A_{\nu^{-i}}^{t_0} \Br \Le (b/\bda)^n,
\en
and hence that 
\eq\label{ADB-LP2}
    \sum_{n\ge1} \sum_{\nu\colon |\nu|=n} \prod_{i=1}^n A_{\nu^{-i}}^{t_0} \ <\ \infty\quad\bbP\mbox{--a.s.}
\en
Furthermore, \Ref{ADB-LP1} also implies that, for all~$n$ large enough, $\bbP(E_n) \le (b/\bda)^n$, where 
\[
     E_n \Def \Blb \sum_{\nu\colon |\nu|=n} \prod_{i=1}^n A_{\nu^{-i}}^{t_0} \ge 1 \Brb\,.
\]
Thus a.s.~only finitely many of the events~$E_n$ occur, and, on~$E_n^c$, since $0 \le t_0 \le 1$,
\eq\label{ADB-LP3}
     \sum_{\nu\colon |\nu|=n} \prod_{i=1}^n A_{\nu^{-i}}^{t_0} 
                \ \ge\ \sum_{\nu\colon |\nu|=n} \prod_{i=1}^n A_{\nu^{-i}}.
\en
\added{\Ref{sucofi} thus follows from \Ref{ADB-LP2} and~\Ref{ADB-LP3}, 
proving~(ii).}
\ep
\end{proofsect}

\section{ A walk that changes its environment, once.}\label{one-change}
\ignore{
Let $A_{\s_{i}}$, with $ i \ge 1$ be as in the previous section. Suppose the kernel of this process satisfies $\ki$. 
}
In this section, we consider a setting in which the process~$\bX$ changes the environment. Fix parameters $L , p > 0$,
and let~$(\added{B}_{\s_{i}},\,i\ge1)$ be a stochastic process, taking values in $[p, +\infty)$, such that the triple
$(A_{\s_{i}}, \added{B}_{\s_{i}}, M_{\s_{i}})$ is a Markov process along rays. 
Recalling that 
$$ 
  S_{\nu} \Def \inf\{ n\ge 0 \colon X_{n} = \nu\},
$$
% where $(X_n,\,n\ge0)$ denotes the random walk run on a ray $[\nu,\infty)$ with states $\nu = \s_0, \s_1,\ldots$,
% For any pair of  adjacent vertices  $\nu$ and  $\nu i$, repeat the following reasoning.  
define 
$$
  G(\nu, n) \Def
    \begin{cases}
     A_{\nu} &\qquad \mbox{ if } \{A_{\nu} > \added{B}_{\nu}\}\cup \{S_{\nu} > n\};\\
     L &\qquad \mbox{ if } \{A_{\nu} \le \added{B}_{\nu}\}\cap \{S_{\nu} \le  n\},
    \end{cases}
$$ 
for each vertex $\nu$ and time $n$.    
If $X_{n} = \nu$, given the environment and $\Fcal_{n} \Def \sigma\{X_{1}, X_{2}, \ldots, X_{n}\}$, 
the probability that  $ X_{n+1} = \nu i$  is given by
\begin{equation}\label{tr1}
    \frac{G(\nu i, n)}{1 + \sum_{j=1}^{b(\nu)}G(\nu j, n)},
\end{equation}
so that the probability of a transition from~$\nu$ to a state~$\nu i$, which has been visited at 
least once before and for which $A_{\nu i} \le \added{B}_{\nu i}$, is modified
by replacing $A_{\nu i}$ by~$L$ in its calculation. As before, the process is absorbed at the state $\parent{\r}$, and recurrence means that the process is absorbed with probability one at $\parent{\r}$.
Let 
$$
  D_{\s_{i}} \Def
    \begin{cases}
      L &\qquad \mbox{ if $A_{\s_{i}} \;\le \; \added{B}_{\s_{i}} $}\\
      A_{\s_{i}} &\qquad \mbox{ if $A_{\s_{i}} \;> \;  \added{B}_{\s_{i}}$.}
    \end{cases}
$$
and denote by $\Kst$ the transition kernel of the Markov chain $\Gamma^{*} \Def 
(D_{\s_{i}}, \added{B}_{\s_{i}}, A_{\s_{i}}, M_{\s_{i}})$ on $\re_+\times\S^*$, where $\S^* := \re_+^2\times\S$
is the state space of $(\added{B}_{\s_{i}}, A_{\s_{i}}, M_{\s_{i}})$, and~$D_{\s_i}$ is singled out. 
\ignore{
Now define $Q^{*}_{C,R}$ as in~\eqref{defnq}, 
with $\Kst$ instead of $K$ and with $\S^*$ playing the part of~$\S$;
set $Q^{*}_{R} := Q^{*}_{-\infty, R}$.
}
As before, define 
$$
    \eta_{\eps,r} \Def 1 - \inf_{y \in \re_+\times\S^*} 
                \mathbb{P}\big(\eps < A_{\s_{1}} \le r \giv  \Gamma_{\s_0}^* = y\big)
$$ 
and $\eta = \lim_{\eps\to 0,r\to\infty} \eta_{\eps,r}$. 

\begin{theorem}\label{rmainthrei}
 Suppose that $\Kst$ satisfies Assumption~$\ki$ and that Assumption~$\kt$ also holds.  Suppose that $L , p \ge 1$.  
Then the condition
\begin{equation}
  \ln (1- \eta)\ >\ -\ln b
           \label{con101} 
\end{equation}
implies the transience of $\mathbf{X}$ on $\Gcal$.
\end{theorem}

\begin{corollary}\label{Cor4.2}
If $\eta = 0$, then the process $\mathbf{X}$ is transient on $\Gcal$.
\end{corollary}

\begin{proofsect}{Proof of Theorem~\ref{rmainthrei}}
Because the process~$\bX$ can change $A_{\nu}$  only at the time~$S_{\nu}$ that~$\nu$ is first visited, 
the proof of Proposition~\ref{prop2} can be used to show that, 
if $\eta<1$ and 
\begin{equation}\label{condmarrei}
 \limsup_{\e \to 0}  \liminf_{n \ti}  \frac1n \int_{\re_+\times\S^*} 
     \ln  \bP\bigl(T_{-1} > T_{n} \;|\; \Gamma^{*}_{\r}= y \bigr) \mmu_{\e}(\d y) 
             \ >\  -\ln b,
\end{equation}
then $\mathbf{X}$ is transient; here, 
$\mmu_{\e}(\cdot)$ is defined as in the previous section, using  $\Kst{}^{\ssup \ell}(x^{*}, \cdot)$ for some $ x^{*} \in \re_+\times \S^*$, and~$\ell$ 
is chosen in such a way that there exist $N $ and~$M$ such that, for all $ x, y \in \re_+\times\S^*$ 
and Borel sets~$B$, we have 
$$ 
    \Kst{}^{\ssup \ell}(x, B) \le \frac MN \sum_{i=1}^{N} \Kst{}^{\ssup i}(y, B).
$$
It remains to determine when \eqref{condmarrei} holds.
 
For a given ray $\s = [\r,\infty)$, let $Q_i^D := \{\sum_{r=0}^{i} \prod_{j=1}^{r} D^{-1}_{\s_j}\}^{-1}$ denote
the probability that the random walk starting in~$\r$ would hit~$\s_i$ before~$\r^{-1}$, if the
probabilities were determined solely by the~$D_{\s_i}$, and, for $i\ge1$, let $q_i^D := Q_i^D / Q_{i-1}^D$ denote
the probability that the same random walk starting in $\s_{i-1}$ hits~$\s_i$ before it hits~$\r^{-1}$.
Then, the probability~$q_i^A$ that the original walk, after it reached $\s_{i-1}$, hits $\s_i$ before~$\r^{-1}$, when 
started in~$\s_{i-1}$, is given by $A_{\s_i}/\{1+A_{\s_i} - q_{i-1}^D\}$, $i\ge1$, with $q_{0}^D$ taken 
to be zero. This leads us to consider the quantity
$$
 \Phi_{n} \Def \prod_{i=1}^n (q_i^A/q_i^D) \Eq 
   \prod_{\heap{i=1}{A_{\s_{i}} < b_{\s_{i}}}}^{n}  
    \left(\frac{1 + D_{\s_i}^{-1}(1 - q_{i-1}^D)}{1 + A_{\s_i}^{-1}(1 - q_{i-1}^D)}\right).
% \prod_{\heap{i=1}{A_{\s_{i}} \in B_{\s_{i}}}}^{n-2}  
%  \left(\frac{\sum_{r=0}^{i-1} \prod_{j=r}^{i-2} D_{\s_j} + L^{-1}}
%   {\sum_{r=0}^{i-1}\prod_{j=r}^{i-2} D_{\s_j} + A^{-1}_{\s_{i}}}\right). 
$$
% We prove that  for any $\eps>0$, there exists a constant $c_{\eps}>0$ such that 
% \begin{equation}
%   \Phi_{n} \1_{\cap_{i=1}^{n}\{A_{\s_{i}}> \eps\}}\ \ge\ c_{\eps} \1_{\cap_{i=1}^{n}\{A_{\s_{i}}> \eps\}}.
%    \label{ADB-Phi-1}
% \end{equation}
Now, since  $D_{\s_i} \ge \theta :=  p \wedge L \ge  1$ for all~$i$, we have 
\eq\label{ADB-q-bnd}
    1 - q_i^D \Le i^{-1},
\en
so that, on the event $\bigcap_{i=1}^{n}\{A_{\s_{i}}> \eps\}$,
\eq\label{ADB-Phi-2}
    \Phi_{n} \ \ge\ \prod_{i=1}^n \{1 +  i^{-1}\eps^{-1}\}^{-1}\ \ge\  k n^{- 1/\e},
\en
for a suitable $k$,  which depends on $\eps$ only, and 
$$ 
    \prod_{i=1}^n q_i^{D} \Eq \prod_{i=1}^n \{1-(1- q_i^{D})\}\ \ge\ 1/n .
$$

Hence, for~\eqref{condmarrei}, we have
\begin{equation}\label{trp}
\begin{aligned}
  \bP(T_{-1} > T_{n} &\giv \Gamma^{*}_{\r} =y) 
    \Eq \mathbb{E}\Bigl[\prod_{i=0}^{n-1} 
          \bP_{\omega}(T_{-1} > T_{i+1} \giv T_{-1} > T_{i}) \Giv \Gamma^{*}_{\r} =y \Bigr]  \\
   &=\ \mathbb{E}\Bigl[\prod_{i=1}^{n} q_i^A \Giv \Gamma^{*}_{\r} =y \Bigr]  
     \Eq \mathbb{E}\Bigl[\Phi_n \prod_{i=1}^{n} q_i^D  \Giv \Gamma^{*}_{\r} =y \Bigr]  \\
   &\ge\ k n^{- 1/\e}\, \mathbb{E}\Bigl[\prod_{i=1}^{n} q_i^D \1_{\cap_{i=1}^{n}\{A_{\s_{i}} >\eps\}}
        \Giv \Gamma^{*}_{\r} =y \Bigr]  \\
   &\ge\ k n^{- (1/\e)-1}   \mathbb{P}\Bigl[\bigcap_{i=1}^{n}\{A_{\s_{i}} >\eps\} \Giv \Gamma^{*}_{\r} =y \Bigr]\,.
\end{aligned}
\end{equation}
Hence, from the definition of~$\eta_{\eps,\infty}$,
$$
   \liminf_{n\ti}\inf_{y \in \re_+ \times \S^*} \frac 1n \ln \bP(T_{-1} > T_{n} \giv \Gamma^{*}_{\r} =y)
   \ \ge\ \ln(1 - \eta_{\eps,\infty}),
$$
and the theorem follows by letting $\e\to 0$ and using~\Ref{condmarrei}.
\hfill\qed
\end{proofsect}

\begin{remark} 
{  Consider a once-reinforced random walk on a Galton--Watson tree, defined as follows. Each edge is 
initially assigned weight 1. The walk moves to any one of its nearest neighbours, with probability 
proportional to the 
weight of the edge traversed. The first time an edge is traversed, its weight becomes $1+\Delta$, for 
$\Delta > -1$, and is never changed again.  With the choice of $L=1$ and $A_{\nu} = 1/(1+\Delta)$ for all
$\nu\in\Gcal$, and with $b_{\nu}= p = \min\{1,1/(1+\Delta)\}$, for all $\nu$, our 
walk is exactly a once-reinforced random walk. Theorem~\ref{rmainthrei} then implies
transience for this class of processes, as already 
proved in \cite{Co06} or \cite{Dai05}. }
\end{remark}
\medskip

{The next result holds for all choices of $L$ and $p$ such that $ L < p$.  } Define 
 the kernel
$$ 
   \Kln^* ((\tw , c, \tu, z),B) \Def K^*((\erm^{\tw}, c, \erm^{\tu},z), E^*(B)),
$$
where $\tw \in [\ln L, \infty)$,  $ c \in (p, \infty)$, $\tu  \in \R,\ z\in\S,\ B \in \Bcal$ and 
$E^*(B) := \{(\erm^{\tw}, c, \erm^{\tu},z)\colon\,(\tw, c, \tu, z)\in B\}$.  Note that, 
if $K^*$ satisfies Assumption~$\ki$, then so does the 
kernel~$\Kln^*$.

{ For  $R  \in (0, \infty)$ and with $C := \ln L$, define the probability kernel $Q^*_{C, R}$ 
on $[C ,R]\times [p, \infty) \times [C ,R] \times \S$ by}
\begin{equation}\label{defnq-star}
    Q^*_{C, R}\big(\ty, (\d \tw, \d c, \d \tu, \d z)\big) 
             \Def \frac{\Kln^* \big(\ty, (\d \tw, \d c,\d \tu, \d z)\big)}
                       {\Kln^*(\ty, [C,  R] \times (p, \infty) \times [C,R] \times \S)}.
\end{equation}
This kernel describes the distribution of the jumps of the process $(\ln D_{\s_i}, \added{B}_{\s_i}, \ln A_{\s_i}, M_{\s_i})$ 
when $\ln A_{\s_i}$ is conditioned to be in the interval $[C ,R]$.  This also implies that $\ln D_{\s_i}$ takes 
values in the same interval. If $K^*$ satisfies Assumption~$\ki$, then so does the 
kernel~$Q^*_{C, R}$. Define 
$$ 
  \tilde{\Lambda}^{\ssup {Q^{*}_{C, R}}} \Def  \limsup_{n \ti} \sup_{\ty } \frac 1n 
          \ln \bbE \left[\erm^{\lambda \sum_{i=1}^n \ln D_{\s_i} } \giv \Gamma^* = \ty  \right],
$$
where the expected value is taken with respect to the kernel  $Q^{*}_{C, R}$, and the supremum over the set
$ [C,  R] \times [p, \infty) \times [C ,R] \times \S$.

{ 
\begin{theorem}\label{rmainthrei2}
 Suppose that $\Kst$ satisfies Assumption~$\ki$ and that  $\eta_{L,\infty} <1$ and $ L < p$.  Then the condition
\begin{equation}
  \limsup_{R \uparrow \infty}\inf_{\lambda \in [0,1]} \tilde{\Lambda}^{\ssup {Q^{*}_{\ln L, R}}}(\lambda)
        \  >\ -\ln b - \ln (1- \eta_{L,\infty})\,    \label{con1010} 
\end{equation}
implies the transience of $\mathbf{X}$ on $\Gcal$.
\end{theorem}

\begin{remark}
Suppose that $\eta = 0$ and $K^*$ satisfies Assumption~$\ki$.  In this case, if $L, p \ge 1$, then, no matter 
what is  the distribution of the  initial environment $(A_{\nu},\,\nu\in\Gcal)$,  the process $\bX$ is transient,
by Corollary~\ref{Cor4.2}.  If instead we 
assume that $L< p$ and $\eta_{L,\infty}<1$, then the process can also be recurrent. In this case, \eqref{con1010} 
provides a sufficient condition for transience.  

\end{remark}
}
%\added{
%\begin{example}
%Consider the binary tree, and let  $A_{\nu}$ be i.i.d. with the following distribution
%$$ A_{\nu} = 
%\begin{cases}
%1/5 \qquad 1/2\\
%3/4 \qquad 1/2.
%\end{cases}
%$$
%\end{example}
%Notice that 
%$$\inf_{\lambda \in [0,1]} \bbE[A^{\lambda}_{\nu}] =\bbE[A_{\nu}] < 1/2.$$
%If we choose $L=p = 1/5$, then the environment does not change. Hence in virtue of Corollary~\ref{impco}, ii), we have that  $\mathbf{X}$ is recurrent. Notice that for this case, $\eta = 0$. Hence, if we choose $L = p =1$, then, in virtue of the previous remark, the random walk $\mathbf{X}$ is transient.
%}
\begin{proofsect}{Proof of Theorem~\ref{rmainthrei2}}  { First, note that 
$$ 
   \Phi_n  \ge \prod_{\heap{i=1}{A_{\s_{i}} \le \added{B}_{\s_{i}}}}^{n}  
    \left(\frac{1 + D_{\s_i}^{-1}(1 - q_{i-1}^D)}{1 + A_{\s_i}^{-1}(1 - q_{i-1}^D)}\right)\1_{\{A_{\s_{i}} \ge  L \}}
   \ \ge\ \prod_{\heap{i=1}{A_{\s_{i}} \le \added{B}_{\s_{i}}}}^{n}  \1_{\{A_{\s_{i}} \ge  L\}},
$$
and hence that
\begin{equation}\label{trp-2}
\begin{aligned}
\bP(T_{-1} < T_{n} \giv \Gamma^{*}_{\r} =y)\ &\ge\  
       \mathbb{E}\Bigl[\prod_{i=1}^{n} q_i^D \1_{\cap_{i=1}^{n-1}\{A_{\s_{i}} \ge  L\}}
        \Giv \Gamma^{*}_{\r} =y \Bigr] \\
&\ge\   \mathbb{E}\Bigl[ \Bigl(\sum_{r=0}^{n}\prod_{j=1}^{r-1} D^{-1}_{\s_n} \Bigr)^{-1} 
        \1_{\cap_{i=1}^{n-1}\{A_{\s_{i}} \ge  L\}}  \giv\, \Gamma^{*}_{\r} =y \Bigr].
\end{aligned}
\end{equation}
}

Now { the last line of \eqref{trp-2} is at most 
$$
\begin{aligned}
    \mathbb{E}\Big[  &\Bigl( (n+1) \max_{r \le n}\prod_{j=1}^{r-1} D^{-1}_{\s_n} \Bigr)^{-1} 
                       \1_{\cap_{i=1}^{n-1}\{A_{\s_{i}} \ge  L\}}  \giv\, \Gamma^{*}_{\r} =y  \Bigr]\\
  &\ge\  \mathbb{E}\Bigl[  \Bigl( (n+1) \max_{r \le n}\prod_{j=1}^{r-1} (D_{\s_n}\wedge R)^{-1} \Bigr)^{-1} 
                       \1_{\cap_{i=1}^{n-1}\{A_{\s_{i}} \ge  L\}}   \giv\, \Gamma^{*}_{\r} =y \Bigr].
\end{aligned}
$$
This, in turn, implies that
$$
\begin{aligned}
 \liminf_{n \ti} \inf_{y \in \S'} \frac 1n \ln &\bbP\bigl(T_{-1} > T_{n} \;|\; \Gamma^{*}_{\r}= y \bigr) \\
    &\ge   \liminf_{n \ti} \inf_{y \in \S'} \frac 1n 
         \ln \bbE\bigl[ {\rm e}^{\min_{t \in (0,1)} \sum_{i=1}^{[nt]} \ln (D_{\s_{i}} \wedge R) }  
              \,\1_{\cap_{i=1}^{n-1}\{A_{\s_{i}} \ge  L\}}\;|\; \Gamma^{*}_{\r}= y \bigr].
\end{aligned}
$$
We now argue much as for~\Ref{new-star} in the proof of the first part of Theorem~\ref{theorem2}, 
% with $D_{\s_{i}}$ instead of $A_{\s_{i}}$, 
proving that 
$$ 
   \liminf_{n \ti} \inf_{y \in \S'} \frac 1n \ln \bbP\bigl(T_{-1} > T_{n} \;|\; \Gamma^{*}_{\r}= y \bigr)
    \ \ge\ \limsup_{R \to \infty}\inf_{\lambda \in [0,1]} \Lambda^{\ssup {Q^{*}_{\ln L,R}}}(\lambda)
              + \ln (1- \eta_L)\ >\  -\ln b.
$$
Hence, \eqref{condmarrei} holds, and this ends the proof.
 \hfill\qed}
\end{proofsect}

\medskip

 As an example, we consider the case where  $\added{B}_{\nu} = p $ is  constant for all $\nu\in\Gcal$.  
Suppose that $L^{-1} = p^{-1} + \eps$ and that $A_{\s_{i}} \in  (L, C)$~a.s., for all $\nu\in\Gcal$  and for
some constant $C$. Note that then $L < p$ and that $\eta_{L,\infty} = 0$.
We prove that~$\bX$ is transient if $\inf_{\lambda \in [0,1]} \Lambda^{\ssup {{K}^{*}_{\ln}}}> - \ln b $, 
and recurrent if $\inf_{\lambda \in [0,1]} {\Lambda}^{\ssup {{K}^{*}_{\ln}}} < - \ln b -\ln (1+ \eps)$.

\added{The transience, when $\inf_{\lambda \in [0,1]} \Lambda^{\ssup {{K}^{*}_{\ln}}}> - \ln b $, }  
is a consequence of Theorem~\ref{rmainthrei2} with \added{$ Q^{*}_{\ln L, \ln C} = K^{*}_{\ln}$}.  
For the proof of recurrence under the assumption 
$\inf_{\lambda \in [0,1]} {\Lambda}^{\ssup {{K}^{*}_{\ln}}} < - \ln b -\ln (1+ \eps)$, 
we have that 
\eqs
  \Phi_{n} &=& \prod_{\heap{i=1}{A_{\s_{i}} \le \added{B}_{\s_{i}}}}^{n} 
     \Bigl\{1 + \frac{(L^{-1} - A^{-1}_{\s_{i}})(1-q_{i-1}^{D})}{1+A^{-1}_{\s_{i}} (1-q_{i-1}^{D})}\Bigr\} \\
   &\le& \prod_{\heap{i=1}{A_{\s_{i}} \le p}}^{n} ( 1+  L^{-1} - A^{-1}_{\s_{i}})\ \le\  
       \prod_{\heap{i=1}{A_{\s_{i}} \le p}}^{n} (1 + \eps) \Le (1+\eps)^{n}.
\ens
Hence
$$
\begin{aligned}
    \bP(T_{-1} > T_{n} &\giv \G^{*} =y) \Eq \mathbb{E}\Bigl[\Phi_n \prod_{i=1}^{n} q_i^D  \Giv \G^{*} =y  \Bigr] 
    \Le (1+\eps)^{n} \cdot \mathbb{E}\Bigl[ \prod_{i=1}^{n} q_i^D  \Giv \G^{*} =y \Bigr].
\end{aligned}
$$

This, { by Theorem~\ref{LDPL}, using Varadhan's lemma and Proposition~\ref{varfor1}, implies that 
$$
    \limsup_{n \ti} \sup_{y \in \S'} \frac 1n \ln \bbP\bigl(T_{-1} > T_{n} \;|\; \G^{*} =y \bigr) 
   \Le \inf_{\lambda \in [0,1]} {\Lambda}^{\ssup {\Kst_{\ln}}}(\l) + \ln (1+ \eps)\ <\ - \ln b.
$$
The expected number of vertices at level $n$ which are visited before the first return to the origin is 
bounded above by $b^{n} \sup_{y \in \S'}  \bbP\bigl(T_{-1} > T_{n} \;|\; \G^{*}= y \bigr)$.
Hence the expected {number of vertices visited before the process returns} to the origin is bounded by
$$ 
   1 + \sum_{n=1}^{\infty } b^{n} \sup_{y \in \S'} \bbP\bigl(T_{-1} > T_{n} \;|\; \G^{*} = y \bigr)\ <\ \infty.
$$
The latter proves  {recurrence}.

\begin{example} \added{With the situation as above, suppose that $(D_{\s_i}, \added{B}_{\s_i}, A_{\s_i})$ 
evolves as a two state Markov chain. $A_\s$ can take the values~$1$ and $p < 1$, and $L^{-1} = p^{-1} + \eps$,
so that $(D_{\s_i}, B_{\s_i}, A_{\s_i})$ has state space $(L, p, p)$ and $(1, p, 1)$.  
We assume that the diagonal elements of the transition matrix  of this process take the value~$3/4$.
Using \eqref{eigtwo}, we have
\[
    \r(\l) \Eq \frac38 \bigl\{ 1 + L^\l + \sqrt{(1+L^\l)^2 - 32L^\l/9} \bigr\},
\]
and  $\Lambda^{\ssup {{K}^{*}_{\ln}}}(\l) = \ln\rho(\l)$.
Since $L < p <1$, this implies that
$$
    \inf_{\l \in [0,1]}\Lambda^{\ssup {{K}^{*}_{\ln}}}(\l) = 
           \ln \Bigl\{ \frac38 \bigl( 1 + L + \sqrt{(1+L)^2 - 32L/9} \bigr) \Bigr\}.
$$
Thus, if $b > 4/3$, the process is always transient.  However, for $1 < b < 4/3$,  the process is recurrent if 
\[
      \frac38 \bigl( 1 + L + \sqrt{(1+L)^2 - 32L/9} \bigr) \ <\ \frac1{b(1+\e)},
\]
and is transient if 
\[
      \frac38 \bigl( 1 + L + \sqrt{(1+L)^2 - 32L/9} \bigr) \ >\ \frac1{b}.
\]   
}
 \end{example}

 \section{Appendix}
\begin{proposition}\label{varfor1}
 Suppose that $\phi \colon \R \to [-\infty, +\infty]$ is a convex function, with $\phi(0)=0$. Then
\begin{equation}\label{vv}
    \sup_{f \in \mathcal{AC}} \Big\{ \min_{t \in [0,1]} f(t) -  
        \int_{0}^{1} \sup_{\lambda} \{f'(u) \lambda - \phi(\lambda)\}\, \d u \Big\} 
            \Eq \inf_{\lambda \in [0,1]} \phi(\lambda).
\end{equation}
\end{proposition}

\begin{proofsect}{Proof}
We first prove that the right-hand side of \eqref{vv} is a lower bound.
Let $\phi$ be finite on $F \subset \R$, and let $t^* \in [0,1]\cap \overline{F}$ be such that
$$
     \lim_{\heap{t \to t^*}{t \in F}}\phi(t) \Eq \inf_{0 \le t \le 1}\phi(t).
$$ 
Such a~$t^{*}$ exists, in virtue of the convexity of $\phi$.
Then, by convexity, $\phi$ has a (non-empty) sub-derivative ${\rm SD}(\phi)\{t^*\}$ at $t^*$. 
Recall that $c \in {\rm SD}(\phi)\{a\}$ means that $\phi(t) \ge \phi(a) + c(t-a)$ for all $t$.

If $t^* \in (0,1)$, then  $ 0 \in SD(\phi)\{t^*\}$, and we choose $f(t) = 0$ for all~$t$ 
to get
$$   
    \inf_{\lambda \in \R} \phi(\lambda) \Eq  \inf_{\lambda \in [0,1]}\phi(\lambda)
$$
as a lower bound for the left hand side of~\eqref{vv}.

If $t^* = 0$, then there is a $c \ge 0$ with $c \in {\rm SD}(\phi)\{0\}$, so that
$\phi(t) \ge \phi(0) + ct$ for all~$ t$.  Take  $f(t) = ct$ for all~$t$. Since $c \ge 0$, 
we have $\min_{0 \le t\le 1}f(t) = 0$, and we get
\[
   - \sup_t \{ct - \phi(t)\}\ \ge\ -\sup_t  \{ct - \phi(0) -ct \} \Eq \phi(0) 
             \Eq \inf_{\lambda \in [0,1]}\phi(\lambda)
\]
as a lower bound for the left hand side of~\eqref{vv}.

If $t^* = 1$, then there is a $c \le 0$ with $c \in {\rm SD}(\phi)\{1\}$, so that
$\phi(t) \ge \phi(1) + c(t-1)$ for all~$t$.  Take $f(t) = ct$. As $c \le 0$, we have   
$\min_{0 \le t\le 1}f(t) = c$, and we get
\[
    c - \sup_t\{ct - \phi(t)\} \ \ge\ c - \sup_t\{ct - \phi(1) -c(t-1)\} \Eq c + \phi(1) - c
            \Eq \inf_{\lambda \in [0,1]}\phi(\lambda)
\]
as a lower bound for the left hand side of~\eqref{vv}.
 
Next we turn to the proof of the upper bound. Fix any $t^{*} \in [0,1]$. Notice that,  
for any $f \in \mathcal{AC}$,  we have  $\min_{t \in [0,1]} f(t) \le 0$ and 
$\big(f(1)- \min_{t \in [0,1]} f(t)\big) \ge 0.$  Hence,
taking $\lambda = t^*$ for all $u\in[0,1]$,
%  $\min_{t \in [0,1]} f(t) < 0$, and 
the left-hand side of  \eqref{vv} is bounded above by 
$$
 \begin{aligned}
   \sup_{f \in \mathcal{AC}} &\big\{ \min_{t \in [0,1]} f(t)  - f(1) t^{*} + \phi(t^{*})\big\}\\
   &\Eq \sup_{f \in \mathcal{AC}} \big\{\min_{t \in [0,1]} f(t) (1- t^{*})  
     - \big(f(1)- \min_{t \in [0,1]} f(t)\big)  t^{*} + \phi(t^{*})\big\} \Le \phi(t^{*}).
 \end{aligned}
$$
By taking the infimum over $ t^{*} \in [0,1]$ we have the upper bound. 
%\qed
\ep
\end{proofsect}

%\begin{theorem}\label{mark1} Suppose $\mathbf{Z} = \{Z_{i}, i \ge 0\}$ is a Markov Chain which takes values in a general state space subset of a Polish space and has transition kernel $K$. 
%Define 
%$$ (K_{\lambda} u)(\sigma) = \sigma^{\lambda} \int_{\Omega} u(\tau) K(\sigma, \d \tau).$$
%Then for all probability measures $\nu$ on $\Omega$, we have  that $(1/N) \sum_{i=1}^{N} \delta_{Z_{i}}$ satisfies LDP with rate function
%$$ \Lambda^{*}(x) \df   \sup_{\lambda} \lambda x -   \Lambda(\lambda),$$
%where 
%$$
% \Lambda(\lambda) \df \lim_{n \ti} \frac 1n \ln \langle (K_{\lambda})^{n} 1, P_{1} \rangle
%$$
%\end{theorem} 

%\begin{theorem} Using the notation of Theorem~\ref{mark1}, we have that 
%$(1/N)\sum_{i=1}^{[Nt]} \ln Z_{i}$ satisfies an LDP with good rate function
%$$ \widetilde{\Lambda}(f) \df \int_{0}^{1} \Lambda^{*}\big(f'(u)\big) \d u.$$
%\end{theorem}

\ignore{
Conceptual Proofs of L Log L Criteria for Mean Behavior of Branching Processes
Russell Lyons, Robin Pemantle, and Yuval Peres
Source: Ann. Probab. Volume 23, Number 3 (1995), 1125-1138. 
}


\begin{thebibliography}{99}

%\bibitem{bill68}
%P. Billingsley (1968), \emph{Convergence of probability measures}, Wiley.

\bibitem{Co06}
A.\ Collevecchio (2006).  
On the transience of processes defined on Galton-Watson trees. 
 {\it Ann.\ Probab\/}. {\bf 34},  870--878.

%\bibitem{BD04}
%B. Davis \& S. Volkov.(2004)  Vertex-reinforced jump processes on trees and finite graphs. {\it Probab. Theory Related Fields} {\bf 128}, Number 1, 42--62. 
\bibitem{Dai05}
J.\ J.\ Dai (2005). 
A once edge-reinforced random walk on a Galton–Watson tree is transient.
{\it Statist.\ Probab.\ Lett.\/} {\bf 73}, 115--124.

\bibitem{dembo-zajic}
A.\ Dembo \& T.\ Zajic (1995). 
Large deviations: from empirical mean and measure to partial sums process, 
{\it Stoch.\ Procs.\ Applics.\/} {\bf 57}, 191--224.

\bibitem{dembo88}
A.\ Dembo \& O. Zeitouni (1998).
 {\it Large deviations techniques and applications\/}, 2nd edition. Springer,
Berlin.

\bibitem{ds89}
J.-D.\ Deuschel \& D.\ W.\ Stroock (1989). 
{\it Large deviations\/}. 
Academic press, Boston.

\bibitem{Du02}
R.\ Durrett, H.\ Kesten \& V.\ Limic (2002). 
Once edge-reinforced random walk on a tree. 
{\it Probab.\ Theory Rel.\ Fields\/}  {\bf 122}, 567--592.

\bibitem{Far11}
G.\ Faraud (2011). 
A central limit theorem for random walk in a random
environment on marked Galton--Watson trees. 
{\it Electronic J.\ Probab.\/} {\bf 16}, 174--215.

\added{
\bibitem{KSK}
J.\ G.\ Kemeny, J.\ L.\ Snell \& A.\ W.\ Knapp (1966). 
{\it Denumerable Markov chains.\/}
van Nostrand, Princeton, NJ.
}
\added{
\bibitem{KS16}
D.\ Kious \&  V.\ Sidoravicius (2016). 
Phase transition for the once-reinforced random walk on $\Z_d$-like trees.  
{\tt arXiv:1604.07631}
}

%\bibitem{ja89}
%P. Jagers.  General branching processes as Markov fields (1989) \emph{Stoch. Proc. Appl.} {\bf 32}, 183 - 212.

\bibitem{LyP92}
R.\ Lyons \& R.\ Pemantle (1992).  
Random walk in a random environment and first--passage percolation on trees . 
{\it  Ann.\ Probab.\/}  {\bf 20}, 125--136.

\bibitem{LyP05}
R.\ Lyons \& Y.\ Peres (2005). 
{\it Probability on trees and networks}. 
Cambridge Series in Statistical and Probabilistic Mathematics
(to appear).


\bibitem{MeP}
M.\ V.\ Menshikov \& D.\ Petritis (2002). 
On random walks in random environment on trees and their relationship with multiplicative chaos. 
{\it Mathematics and computer science II (Versailles, 2002)\/},
 415--422. 

\bibitem{nn1987}
P.\ Ney \& E.\ Nummelin (1987). 
Markov additive processes I: Eigenvalue properties and limit theorems. 
{\it Ann.\ Probab.\/} {\bf 15}, 561--592.

\bibitem{nn1987b)}
P.\ Ney \& E.\ Nummelin (1987). 
Markov additive processes II: Large deviations. 
{\it Ann.\ Probab.\/} {\bf 15}, 593--609.

\bibitem{Pe07}
R.\ Pemantle (2007). 
A survey of random processes with reinforcement. 
{\it Probab.\ Surv.\/} {\bf 4}, 1--79.
 
\bibitem{SZ02}
A.\ Sznitman (2002). 
Topics in random walks in random environment. 
In {\it School and Conference in Probability Theory\/}, ICTP Lect. Notes, XVII 203--266 (electronic). 
Abdus Salam Int.\ Cent.\ Theoret.\ Phys., Trieste. 
Available at {\tt http://users.ictp.trieste.it/~pub\_off/lectures/lns017/Sznitman/Sznitman.pdf}

\end{thebibliography}
\end{document}